\documentclass[preprint,12pt]{elsarticle}

\usepackage{amsmath,amsfonts,amssymb,dsfont,mathrsfs}
\usepackage{amssymb}
\usepackage{latexsym,stmaryrd}
\usepackage{pstricks}
\usepackage{textcomp}
\usepackage{psfrag}
\usepackage{graphicx}

\def\Ree{\mbox{$\mathds{R}$}}

\def\Cee{\mbox{$\mathds{C}$}}
\def\Nee{\mbox{$\mathds{N}$}}
\def\Zee{\mbox{$\mathds{Z}$}}

\usepackage{hyperref}

\journal{Systems \& Control Letters}

\begin{document}

\begin{frontmatter}
\title{Stability analysis for a class of linear systems governed by difference equations}
\author{S\'erine Damak\corref{CorresAuthor}\fnref{CSCship}}
\ead{serine.damak@insa-lyon.fr}
\author{Michael Di Loreto\fnref{CSCship}}
\ead{michael.di-loreto@insa-lyon.fr}
\address{Laboratoire Amp\`ere, UMR CNRS 5005, INSA-Lyon, 20 Avenue Albert Einstein, 69621 Villeurbanne, France}

\author{Warody Lombardi}
\ead{warody.lombardi@cea.fr}
\address{CEA-LETI, Minatec Campus, 17 rue des Martyrs, 38054 Grenoble Cedex, France.}

\author{Vincent Andrieu}
\ead{vandrieu@lagep.univ-lyon1.fr}
\address{Universit\'e de Lyon, LAGEP, 43 Bd du 11 novembre 1918, 69621 Villeurbanne, France.}

\cortext[CorresAuthor]{Corresponding author}
\fntext[CSCship]{The authors acknowledge the financial support of the French National Research Agency under ANR project entitled Approximation of Infinite Dimensional Systems.}

\begin{abstract}
Linear systems governed by continuous-time difference equations cover a 
wide class of linear systems. From the Lyapunov-Krasovskii approach, we investigate stability for such a class of systems. Sufficient conditions, and in some particular cases, necessary and sufficient conditions for exponential stability are established, for multivariable systems with commensurate or rationally independent delays. A discussion on robust stability is proposed, for parametric uncertainties and time-varying delays. 
\end{abstract}

\end{frontmatter}

\section{Introduction}\label{section1}
In this note, we are interested with the class of linear systems governed by continuous-time difference equations described by  
\begin{equation}\label{eq1-1}
x(t)= \sum_{k=1}^{N}A_k \, x(t-r_k),
\end{equation} 
where $x(t)\in \Ree^{n}$ is called the instantaneous state at time $t\geq 0$, $A_{k}$ are real $n\, \times \, n$ matrices, for $k=1,\ldots,N$, and ($r_1,\ldots, r_N$) are the delays, with $0<r_1< r_2< \ldots < r_N$.\\
The motivations to work on such a class of systems come from conservation laws, neutral time-delay systems, sampled-data systems, or from approximation of distributed-delays.
The system~(\ref{eq1-1}) is a particular case of the renewal equation presented in~\cite{bellmancooke}, where some general conditions for the existence and the unicity of a solution were established. After this preliminary work, stability analysis for~(\ref{eq1-1}) was a central topic of many researches, with a particular emphasis on robust stability for small variations in the delays. Based on a functional analysis, spectral conditions for stability independent of the delays were proposed in~\cite{melvin} for the scalar case, and in~\cite{silkowskii} or~\cite{avellarhale} for the multivariable case. Variations in the delays were also studied for state feedback control, as in~\cite{haleverduynlunelSmallDelays},~\cite{haleverduynlunel2003}. Stability conditions were also obtained from Lyapunov-Krasovkii techniques. In~\cite{carvalho}, the author investigated stability and asymptotic stability for~(\ref{eq1-1}), handling out conditions expressed in terms of Linear Matrix Inequalities (LMI). A construction of Lyapunov functionals for~(\ref{eq1-1}) was proposed in~\cite{shaikhet}. The second method of Lyapunov was analyzed in~\cite{pepe} for a more general class of nonlinear difference equations.\\
Stability for various extensions of~(\ref{eq1-1}) was also studied. A first contribution was proposed in~\cite{avellarmarconato}, where a sufficient condition for asymptotic stability for time-varying parameters and delays in~(\ref{eq1-1}) was established in the scalar case. The authors outlined that such an extension in the multivariable case was not trivial. A first answer on stability for time-varying delays was positively discussed in~\cite{damak2013}. Other extensions of classes of systems include the works on neutral time-delay systems (see for instance~\cite{halelunelbook},~\cite{fridman},~\cite{kharitonov2005},~\cite{gukharitonovchen} or~\cite{fridmanLPV}), or systems with distributed delays~\cite{mondieAguilar},~\cite{melchorAguilar}. \\
From the Lyapunov-Krasovskii approach, we propose in this paper some new sufficient conditions for exponential stability of~(\ref{eq1-1}). For the case of commensurate delays, necessary and sufficient conditions for exponential stability are characterized. These conditions are LMI, for which efficient numerical algorithms exist. Estimations of the exponential decay rate are proposed, allowing to extend the purpose of the conditions given in~\cite{carvalho}. A discussion on robustness under parametric norm-bounded uncertainties is made. The last contributions include sufficient conditions for exponential stability for time-varying delays.\\
The paper is organized as follows. Section~\ref{section2} addresses some properties of the solution for~(\ref{eq1-1}), and basic definitions are recalled. In Section~\ref{section3}, we briefly present some known results on stability. Exponential stability is solved in~Section~\ref{section4}. Robustness for parametric uncertainties is discussed in Section~\ref{section5}, while Section~\ref{section6} addresses stability for time-varying delays. Examples with simulations illustrate the various conditions on stability. 

Let us introduce few notations. For any bounded continuous initial function $\varphi \in \mathcal{C} ([-r_N,0[,\Ree^n)$, the solution $x(t,\varphi)$ of~(\ref{eq1-1}) with initial condition $\varphi(\cdot)$ is well defined and unique for $t\geq 0$~\cite{bellmancooke}. Such a solution will be denoted sometimes by $x(t)$, if no confusion on the initial condition dependency arises. We denote by $x_t(\varphi)$ the partial state trajectory, for $t\geq 0$, that is
$$
x_t(\varphi):\theta \mapsto x(t+\theta,\varphi)\; , \; \theta\in[-r_N,0[.
$$
The space of initial continuous functions is endowed with the norm $|| \varphi||_{c}=\max_{\theta \in [-r_{N},0]} ||\varphi(\theta) ||$, where $||\cdot||$ stands for the Euclidean norm. $\|x_t(\varphi)\|_{L_2}$ stands for the $L_2$-norm, that is
$$
\|x_t(\varphi)\|_{L_2}^2 = \int_{-r_N}^{0}{\|x(t+\theta,\varphi)\|^2\,\mathrm{d}\theta}.
$$
We denote by $P^T$ the transposed matrix of $P$, and $\lambda_{\min}(P)$ (resp. $\lambda_{\max}(P)$) the smallest (resp. the largest) eigenvalue of a  symmetric positive definite matrix $P$, that we will abbreviate by $P>0$, or by $P\geq 0$ if the matrix $P$ is positive semidefinite. Similar notations will be used for symmetric negative definite (resp. semidefinite) matrices. The spectral radius of a matrix $A$ is denoted by $\rho(A)$.

\section{Systems governed by difference equations}\label{section2}
\subsection{Properties on discontinuities}\label{section2-1}
For any bounded continuous initial function $\varphi \in \mathcal{C} ([-r_N,0[,\Ree^n)$, the solution $x(t,\varphi)$ of~(\ref{eq1-1}) is piecewise continuous, for all $t\geq 0$. The computation of this solution is obtained through a direct time-recursive scheme, which reproduces linear combinations of the past solution in time. The discontinuities are propagated from time $t_0=0$, where
$$
x(t_0) = \sum_{k=1}^{N}{A_k\,\varphi(-r_k)}
$$
is, in general, different from $\varphi(t_0^-)=\mathrm{lim}_{t\rightarrow 0^-}\varphi(t)$. We will denote $t_k$, for $k\in\Nee$ the times of these discontinuities. For any $k\in\Nee$, the solution $x(t,\varphi)$ is continuous over $[t_k,t_{k+1}[$. By iterating in time the solution, it is readily verified that the time discontinuities $t_k$ are governed by the recursive formula
$$
t_{k} = \mathop{\mathrm{min}}_{m_{k}^1,\ldots,m_{k}^N}{\left\{\sum_{i=1}^{N}{m_{k}^ir_i}\; : \; t_{k}>t_{k-1}, \; m_{k}^i\in\Nee\right\}},
$$
with $t_0=0$. The time $\delta_k=t_{k}-t_{k-1}$ between two jump discontinuities in $t_{k-1}$ and $t_{k}$ satisfies 
$$
\delta_k = \mathop{\mathrm{min}}_{m_{k}^1,\ldots,m_{k}^N}{\left\{\sum_{i=1}^{N}{(m_{k}^i-m_{{k-1}}^i)r_i}>0\; , \; m_{k}^i\in\Nee\right\}}.
$$
The delays are said to be rationally independent if $\sum_{k=1}^{N}{m_kr_k}=0$ for some $(m_1,\ldots,m_N)\in\Zee^N$ implies that $m_i=0$ for $i=1,\ldots,N$. When the delays are rationally independent, one can verify that $\frac{r_i}{r_j}$ are irrational numbers, for any $i\neq j$. It is then a direct consequence of Dirichlet theorem to see that 
$\mathop{\mathrm{inf}}_{k\in\mathds{N}}{\,\delta_k} = 0$. But this infimum bound can not be reached, by definition. This fact is to compare with the case of commensurate delays, that is $r_k=kr$ for some $k\in\Nee$ and $r>0$,  for which $\delta_k = r$ for any $k\in\Nee$. In this last case, the successive jump discontinuities arise at times $t_k = kr$, for $k\in\Nee$.
\subsection{Stability}\label{section2-2}
From these basic remarks on the discontinuity of the solution, let us recall some definitions of stability for systems in the form~(\ref{eq1-1}).

\newdefinition{defn}{Definition}
\begin{defn}\label{def2-2-1}
The system~(\ref{eq1-1}) is said to be
\begin{enumerate}
\item[$i)$] stable (resp. $L_2$-stable) if, for any $\epsilon>0$, there exists $\delta(\epsilon)>0$ such that 
$\|\varphi\|_c < \delta$ implies that $\|x(t,\varphi)\| < \epsilon$ (resp. $\|x_t(\varphi)\|_{L_2} < \epsilon$), for any $t \geq 0$.
\item[$ii)$] $L_2$-asymptotically stable if it is $L_2$-stable, and for any bounded initial function $\varphi$ in $\mathcal{C}([-r_N,0[,\Ree^n)$, 
$$
\underset{t \rightarrow \infty}{\lim} \, \|x_t(\varphi)\|_{L_2} = 0.
$$
\item[$iii)$] asymptotically stable if it is stable, and for any bounded initial function $\varphi$ in $\mathcal{C}([-r_N,0[,\Ree^n)$, 
$$
\underset{t \rightarrow \infty}{\lim} \, x(t,\varphi) =0.
$$
\item[$iv)$] $L_2$-exponentially stable if it is $L_2$-asymptotically stable,
and if there exist $\alpha\geq 0$ and $\mu >0$  such that
$$ || x_t(\varphi)||_{L_2} \leq \alpha \,\mathrm{e}^{-\mu t} || \varphi||_c, \; \; \forall \, t \geq 0.$$
\item[$v)$] exponentially stable if it is asymptotically stable,
and if there exist $\alpha \geq 0$ and $\mu >0 $ such that 
$$ || x(t,\varphi)|| \leq \alpha \,\mathrm{e}^{-\mu t} || \varphi ||_c, \; \; \forall \, t \geq 0.$$
\end{enumerate}
\end{defn}
It is clear that exponential stability implies $L_2$-exponential stability, as well as asymptotic stability. However, the converse is false, in general. Furthermore, these definitions are done for a given set of delays $\{r_1,\ldots,r_N\}$. If these properties of stability hold independently of the delays, we will say that~(\ref{eq1-1}) is stable (asymptotically, exponentially) in the delays.

\section{Asymptotic stability analysis}\label{section3}

\subsection{Spectral analysis}\label{section3-1}

Some results are available in the literature in which $L_2$-asymptotic stability is studied for system~(\ref{eq1-1}). In this section, we remind the reader these results. In~\cite{avellarhale} and~\cite{silkowskii}, the authors give a necessary and sufficient condition for $L_2$-asymptotic stability in the delays, that is
\begin{equation}\label{eq3-1-1}
\mathrm{sup}\left\{\rho(\sum_{k=1}^{N}{\mathrm{e}^{j\theta_k}A_k}),\; \theta_k\in[0,2\pi]\right\}<1.
\end{equation}
In the scalar case, this condition comes down to
\begin{equation}\label{eq3-1-2}
\sum_{k=1}^{N}{|A_k|}<1.
\end{equation}
For commensurate delays, the system~(\ref{eq1-1}) admits a state-space realization with a single delay $r>0$ of the form
\begin{equation}\label{eq3-1-3}
x(t) = A \, x(t-r).
\end{equation}
In this particular case, a complete equivalence on stability for linear sampled-data systems holds. See for instance~\cite{carvalho} and~\cite{damak2013}. This stability is of course in the delay $r$, since these conditions are independent of the delay.
\newdefinition{thm}{Theorem}
\begin{thm}\label{th3-1-1} 
The system~(\ref{eq3-1-3}) is 
\begin{itemize}
\item[$i)$] asymptotically stable if and only if $\rho(A)<1$.
\item[$ii)$] stable if and only if $\rho(A)\leq 1$, and for any unit eigenvalue $|\lambda_{k}| =1$, $\mathrm{rank}(A-\lambda_{k}I)=n-q_{k}$, where $q_{k}$ is the algebraic multiplicity of $\lambda_{k}$.
\end{itemize} 
\end{thm}

\newdefinition{exmp}{Example}
\begin{exmp}\label{ex1}
Consider the system
\begin{equation}\label{eq-ex1}
x(t)=\frac{3}{4}x(t-r_1)-\frac{3}{4}x(t-r_2)
\end{equation}
with $r_1=1$, $r_2=2$, and the initial condition $\varphi(t)=2\,\sin(t) $, for $t\in [-2,0[$. This system has commensurate delays, and can be written in the form
$$
X(t) = A\,X(t-r_1),
$$
with $A=\begin{bmatrix}
\frac{3}{4}  &  -\frac{3}{4}\\
1   &  0
\end{bmatrix}$, $X(t)=\begin{bmatrix}
x(t) \\ x(t-r_1)
\end{bmatrix}$.
From Theorem~\ref{th3-1-1}, it is asymptotically stable, since $\rho(A) = \frac{\sqrt{3}}{2}$. However, it is not stable in the delays, since~(\ref{eq3-1-2}) gives
$\frac{3}{4}+\frac{3}{4}=\frac{3}{2}>1$. According to~\cite{avellarhale}, unstability appears for some arbitrarily small variations in the delays. For instance, for $r_2 = 2+\frac{\pi}{10}$, the system~(\ref{eq-ex1}) is unstable, as shown in Fig.~\ref{fig-1} where a simulation result for $x(t,\varphi)$ is provided. 
\begin{figure}[!h]
\begin{center}\psfrag{t}[c]{\scriptsize $t$}
\psfrag{y}[c]{\scriptsize $x(t,\varphi)$}
\includegraphics[width=8cm]{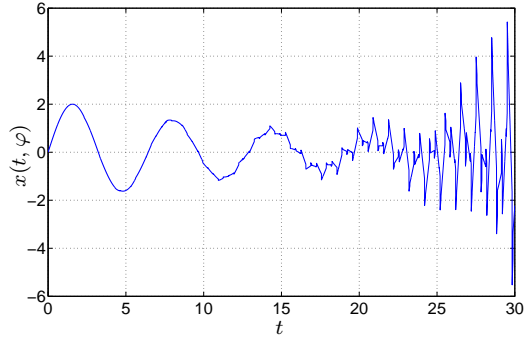}  
\caption{\it Unstability of system~(\ref{eq-ex1}) with $r_1=1$ and $r_2=2+\frac{\pi}{10}$.}\label{fig-1}
\end{center}\end{figure}

\end{exmp}

\subsection{Analysis based on Lyapunov-Krasovskii functionals}\label{section3-2}

For system~(\ref{eq3-1-3}), a necessary and sufficient condition for asymptotic stability can also be obtained by Lyapunov-Krasovskii techniques. Let us synthesize this condition in the following numerical condition.
\begin{thm}\label{th3-2-1} \cite{carvalho}
The system~(\ref{eq3-1-3}) is asymptotically stable if and only if for any given symmetric positive definite real matrix $M$, there exists a symmetric positive definite real matrix $P$ such that
\begin{equation}\label{eq3-2-1}
A^{T}P A -P= -M.
\end{equation}
If~(\ref{eq3-2-1}) is satisfied by a symmetric positive definite matrix $P$ and a symmetric positive semidefinite matrix $M$, then~(\ref{eq3-1-3}) is stable.
\end{thm}

A similar stability condition holds for arbitrary delays $r_k$, $k=1,\ldots,N$. However, this condition is not necessary, and is related to $L_2$-stability.

\begin{thm}\label{th5} \cite{carvalho}
The system~(\ref{eq1-1}) is $L_2$-asymptotically stable if for any given symmetric positive definite real matrix $M$, there exist symmetric positive definite real matrices $P_k$, $k=1,\ldots,N$, such that~(\ref{eq3-2-2}) is fulfilled.\\
\begin{figure*}
\begin{equation}\label{eq3-2-2}
-M = \begin{bmatrix}
A^{T}_{1} P_{1} A_{1} -P_1 + P_{2}     & A_{1}^{T} P_{1} A_{2}   &  \cdots &   A_{1}^{T} P_{1} A_{N}\\
A^{T}_{2} P_{1} A_{1}    &   A_{2}^{T} P_{1} A_{2} + P_{3} - P_2   &   \cdots   &  A^{T}_{2} P_{1} A_{N}  \\
\vdots &  \vdots &  \cdots    &\vdots\\
\vdots & \vdots    &    \ddots    &\\
A^{T}_{N} P_{1} A_{1}  &   A_{N}^{T} P_{1} A_{2} & \cdots   & A^{T}_{N} P_{1} A_{N}-P_N 
\end{bmatrix}.
\end{equation}
\hrulefill
\end{figure*}
If (\ref{eq3-2-2}) is satisfied by symmetric positive definite matrices $P_k$, $k=1,\ldots,N$, and a symmetric positive semidefinite matrix $M$, then~(\ref{eq1-1}) is $L_2$-stable.
\end{thm}

Obviously, the stability characterized in Theorem~\ref{th5} is in the delays. In~(\ref{eq3-2-2}) when $M$ is positive definite, the computation of the exponential decay rate can not be retrieved directly from the previous theorem. Furthermore, more information about the type of stability may be obtained. For this, inspired from~\cite{mondiekharitonov}, we adapt in the next section the Lyapunov-Krasovskii approach to test exponential stability and to compute an exponential decay rate.

\section{Exponential stability}\label{section4}

\subsection{Single-delay case }\label{section4-1}

Let us start with the single-delay case. Consider the system
\begin{equation}\label{eq4-1-1}
x(t)=A \, x(t-r) ,
\end{equation}
where $A\in \Ree^{n\times n}$, $r>0$, and an initial bounded condition defined by $ \varphi \in \mathcal{C}([-r,0[,\Ree^n)$. We have the following stability condition.
We have this following result. 

\begin{thm}\label{th4-1-1}
If there exist a real $n\times n$ symmetric positive definite matrix $P$ and a positive constant $\mu>0$ such that the inequality  
\begin{equation}\label{eq4-1-2}
-M_{\mu}=A^T P A -{\mathrm{e}^{-2\mu r}}P\leq 0
\end{equation}
holds, then~(\ref{eq4-1-1}) is exponentially stable, that is 
\begin{equation}\label{eq3-4-3}
\|x(t,\varphi)\|\leq \sqrt{\frac{\lambda_{\max}(P)}{\lambda_{\min}(P)}}\mathrm{e}^{-\mu t}\|\varphi\|_c \; , \; t\geq 0.
\end{equation}
Conversely, if the system~(\ref{eq4-1-1}) is exponentially stable, there exist $P>0$ and $\mu>0$ such that~(\ref{eq4-1-2}) is satisfied.
\end{thm}
\newproof{pf}{Proof}
\begin{pf}
Consider the Lyapunov-Krasovskii functional
\begin{equation}\label{eq4-1-3}
v_{\mu}(x_t(\varphi))=\int_{t-r}^{t} \mathrm{e}^{-2\mu (t-\theta)} x^{T}(\theta) P x(\theta) \, \mathrm{d}\theta , 
\end{equation}
with $\mu>0$. This functional satisfies 
\begin{eqnarray}
 \lambda_{\min}(P) \mathrm{e}^{-2\mu r}  || x_t(\varphi)||^2_{L_{2}} & \leq & v_{\mu}(x_t(\varphi)) \nonumber,\\
 \lambda_{\max}(P) || x_t(\varphi)||^2_{L_{2}} & \geq & v_{\mu}(x_t(\varphi)).\label{eq4-1-4}
\end{eqnarray} 
The time derivative of $v_\mu(x_t(\varphi))$ along the trajectories of~(\ref{eq4-1-1}) is
\begin{eqnarray*}
 \frac{\mathrm{d}}{\mathrm{d}t}v_\mu(x_t(\varphi)) & = & -2\mu \,v_\mu(x_t(\varphi))- x^{T}(t-r) M_\mu x(t-r).
\end{eqnarray*}
If~(\ref{eq4-1-2}) is satisfied, we conclude that
$$
\frac{\mathrm{d}}{\mathrm{d}t}v_\mu(x_t(\varphi)) +2\mu \, v_\mu(x_t(\varphi))\leq 0 \,,\; \forall \, t\geq 0.
$$
This inequality implies that $v_\mu(x_t(\varphi))\leq \mathrm{e}^{-2\mu t} v_\mu(\varphi)$ for $t\geq 0$. From (\ref{eq4-1-4}), we obtain 
\begin{equation}\label{eq3-4-6}
v_\mu(x_t(\varphi))\leq r \lambda_{\mathrm{max}}(P)\,\|\varphi\|_c^2\,\mathrm{e}^{-2\mu t},\; t\geq 0.
\end{equation}
To conclude on exponential stability, an upper bound for the euclidean norm of $x(t,\varphi)$ need to be established. For this, note that $M_\mu\geq 0$ in~(\ref{eq4-1-2}) implies
$$
x^T(t)Px(t) \leq \mathrm{e}^{-2\mu r}x^T(t-r)Px(t-r),\; t\geq 0.
$$
Iterating this last inequality leads to \cite{kharitonov2005}
$$
0\leq \|x(t,\varphi)\| \leq \alpha\, \mathrm{e}^{-\mu t}\, , \; t\geq 0, 
$$
where
$$
\alpha =\sqrt{ \frac{{\lambda_{\mathrm{max}}(P)}}{{\lambda_{\mathrm{min}}(P)}}}\|\varphi\|_c.
$$
Conversely, if~(\ref{eq4-1-1}) is exponentially stable, it follows from Theorem~\ref{th3-2-1} that~(\ref{eq3-2-1}) is fulfilled, for some positive definite matrices $P$ and $M$. Since $1-\frac{\lambda_{\mathrm{min}}(M)}{\lambda_{\mathrm{max}}(P)}$ is in $[0,1[$, take any $\mu$ in the interval 
$$
0 < \mu \leq -\frac{1}{2r}\,\mathrm{ln}\left(1-\frac{\lambda_{\mathrm{min}}(M)}{\lambda_{\mathrm{max}}(P)}\right).
$$
It is then a routine to verify that, for such a $\mu$, 
$$
-M_\mu = -M+(1-\mathrm{e}^{-2\mu r}) P \leq 0,
$$ 
so that~(\ref{eq4-1-2}) holds.
\hfill{$\Box$}
\end{pf}

\begin{exmp}\label{ex-2}
Let us consider the system
\begin{equation}\label{eq-ex2}
x(t)= \begin{pmatrix}
\frac{1}{2} & -\frac{3}{10} \\ \frac{7}{20} & 0\\
\end{pmatrix}
 x(t-\pi), 
\end{equation}
with initial condition $\varphi(t)=\begin{bmatrix}\sin(3t) \\ \cos(3t)\end{bmatrix}$, for $t\in[-\pi, 0[$. The system~(\ref{eq-ex2}) is exponentially stable, i.e., a solution of~(\ref{eq4-1-2}) is
$$
P = \begin{pmatrix} 22.8565 & -16.3276 \\ -16.3276 & 19.5955\end{pmatrix},\; \mu = 0.3584.
$$
It follows that
$$
\|x(t,\varphi)\| \leq  2.7951\cdot\mathrm{e}^{- 0.3584\, t},\;t\geq 0.
$$
The free response of~(\ref{eq-ex2}) is plotted in Fig.~\ref{fig-2}. 

\begin{figure}[h!]
\center
\psfrag{y}[c]{\scriptsize $x(t,\varphi)$}
\psfrag{t}{\scriptsize $t$}
\includegraphics[width=9cm]{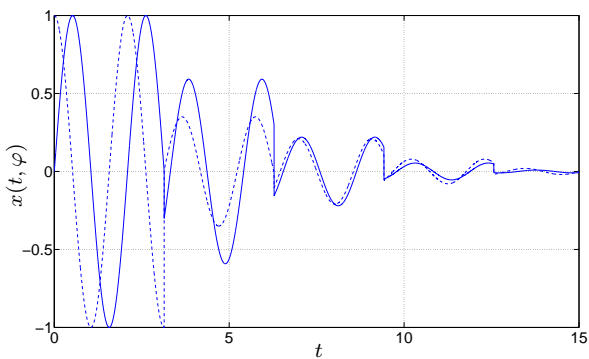}
\caption{Free response $x(t,\varphi)$ of the system~(\ref{eq-ex2}): $x_1(t,\varphi)$ (continuous line), $x_2(t,\varphi)$ (dashed line).}\label{fig-2}
\end{figure}
\end{exmp}

\subsection{Multi-delays case }\label{section4-2}

Let us consider now the case of systems with arbitrary delays
\begin{equation}\label{eq4-2-1}
x(t)=\sum_{k=1}^{N} A_{k} x(t-r_{k}),
\end{equation}
with initial condition $\varphi \in \mathcal{C}([-r_N,0[,\Ree^{n})$. For the multi-delays case, similar arguments than the single-delay case partially hold. Since it is of independent interest, we mention the following general result for stability of~(\ref{eq4-2-1}). The partial trajectory $x_t(\varphi)$ lies in the space $\mathcal{P_C}([-r_N,0[,\Ree^n)$ of piecewise continuous vector functions.

\begin{thm}\label{th4-2-1}
Assume that there exists a continuous functional $v\;:\mathcal{P_C}([-r_N,0[,\Ree^n)\rightarrow \Ree$ such that $t\mapsto v(x_t(\varphi))$ is (upper right-hand) differentiable for all $t\geq 0$ and such that
\begin{enumerate}
\item $\exists\, \alpha_1>0$ s.t. $\forall t\geq 0$, $\alpha_1 \|x_t(\varphi)\|_{L_2}^{2} \leq v(x_t(\varphi))$,
\item $\exists\, \alpha_2\geq 0$ s.t. $v(\varphi) \leq \alpha_2 \|\varphi\|_{c}^{2}$, 
\item $\exists\, \mu>0$ s.t. $\forall t\geq 0$, $\frac{\mathrm{d}}{\mathrm{d}t}v(x_t(\varphi))+2\mu \,v(x_t(\varphi))\leq 0$.
\end{enumerate}
Then~(\ref{eq4-2-1}) is $L_2$-exponentially stable, that is
$$
\|x_t(\varphi)\|_{L_2} \leq \sqrt{\frac{\alpha_2}{\alpha_1}} \|\varphi\|_c\,\mathrm{e}^{-\mu t},\; t\geq 0.
$$
\end{thm}
\begin{pf}
Assumption~(3) leads to
$$
v(x_t(\varphi)) \leq v(\varphi)\,\mathrm{e}^{-2\mu t},\;t\geq 0.
$$
From assumptions~(1) and~(2), we obtain
$$
\|x_t(\varphi)\|_{L_2}^2 \leq \frac{\alpha_2}{\alpha_1} \|\varphi\|_c^2\,\mathrm{e}^{-2\mu t},\;t\geq 0,
$$
which proves the assertion.\hfill{$\Box$}
\end{pf}

The conditions of this result are strongly similar to Theorem 3 in~\cite{melchorAguilar}, but the conclusion is a bit different, since here only $L_2$-exponential stability can be proved through such a result. Applying Theorem~\ref{th4-2-1} for~(\ref{eq4-2-1}) leads to the following sufficient condition for stability.

\begin{thm}\label{th4-2-2}
If there exist symmetric positive definite real matrices $P_k$, $k=1,\ldots,N$, and $\mu>0$, such that $M_{\mu}$ in~(\ref{eq4-2-2}) is a positive semidefinite matrix, then the system~(\ref{eq4-2-1}) is $L_2$-exponentially stable, that is 
$$ 
\|x_t(\varphi) \|_{L_{2}} \leq \alpha\|\varphi\|_{c} \,\mathrm{e}^{-\mu t}  ,\hspace{0.4 cm} \forall t \geq 0,
$$
where 
$$\alpha^2=\frac{\sum_{k=1}^{N} (r_k-r_{k-1})\lambda_{\max}(P_k) }{ \lambda_{\min}(P_N)\,\mathrm{e}^{-2\mu r_N}}. $$
\begin{figure*}
\begin{equation}\label{eq4-2-2}
\small -M_{\mu} = \begin{bmatrix}
A^{T}_{1} P_{1} A_{1} + {\mathrm{e}^{-2 \mu r_{1}}} (P_2 - P_{1})     & A_{1}^{T} P_{1} A_{2}   &   & \cdots  & A_{1}^{T} P_{1} A_{N}\\
A^{T}_{2} P_{1} A_{1}    &   A_{2}^{T} P_{1} A_{2} + {\mathrm{e}^{-2 \mu r_{2}}} (P_{3} - P_2)   &     & \cdots  & A^{T}_{2} P_{1} A_{N}  \\
 \vdots & \vdots  &   & \cdots & \vdots\\
  \vdots   & \vdots &  & \ddots   &\\
A^{T}_{N} P_{1} A_{1}  &   A_{N}^{T} P_{1} A_{2} & \cdots &  & A^{T}_{N} P_{1} A_{N}-{\mathrm{e}^{-2 \mu r_N}} P_N 
\end{bmatrix}.
\end{equation}
\hrulefill
\end{figure*}
\end{thm}
\begin{pf} Assume that there exist symmetric positive definite matrices $P_k$, for $k=1,\ldots,N$, and $\mu>0$ such that $M_\mu$ in~(\ref{eq4-2-2}) is positive semidefinite. Consider the Lyapunov-Krasovskii functional
$$
v_{\mu}(x_t(\varphi))=\sum_{k=1}^{N} \int_{t-r_{k}}^{t-r_{k-1}} \mathrm{e}^{-2\mu (t-\theta)} x^{T}(\theta) P_k x(\theta) \,\mathrm{d}\theta.
$$
The functional $v_\mu(x_t(\varphi))$ is continuous, differentiable with respect to $t$, and satisfies
$$
\alpha_1 \|x_t(\varphi)\|_{L_2}^2 \leq v_\mu(x_t(\varphi)),
$$
for
$\alpha_1 = \min\limits_{k=1,\ldots,N}(\lambda_{\min}(P_k)\,\mathrm{e}^{-2\mu r_k})>0$, and
$$
v_\mu(\varphi)\leq \alpha_2 \|\varphi\|_c^2,
$$
with $\alpha_2 = \sum_{k=1}^{N} (r_k-r_{k-1})\lambda_{\max}(P_k)$. From the assumption on $M_\mu$ in~(\ref{eq4-2-2}), it is straightforward to verify that $P_N\leq P_{N-1}\leq \ldots \leq P_1$. This in turn implies that $\alpha_1 = \lambda_{\min}(P_N)\,\mathrm{e}^{-2\mu r_N}$. Furthermore, its time derivative along the trajectories of~(\ref{eq4-2-1}) is given by
\begin{equation}\label{eq4-2-3}
\frac{\mathrm{d}}{\mathrm{d}t}v_{\mu}(x_t(\varphi))= -2\mu \,v_{\mu }(x_t(\varphi))-\psi^{T}(t) M_{\mu } \psi(t),
\end{equation}
where 
$\psi^{T}(t)=\begin{bmatrix}
x^{T}(t-r_1) & \cdots  & x^{T}(t-r_N)   
\end{bmatrix}$.
Then 
$$
\frac{\mathrm{d}}{\mathrm{d}t}v_\mu(x_t(\varphi)) +2\mu\, v_\mu(x_t(\varphi))\leq 0\, , \; t\geq 0. 
$$
The result follows from Theorem~\ref{th4-2-1}.\hfill{$\Box$}
\end{pf}
Similarly to Corollary 5.4 in~\cite{carvalho}, we have the following corollary, where the link between the conditions on exponential stability in the delays appears.

\newdefinition{cor}{Corollary}
\begin{cor}\label{th4-2-3}
For the system~(\ref{eq4-2-1}), if there exist positive definite symmetric matrices $P_k>0$, for $k=1,\ldots,N$ and $\mu>0$ such that the symmetric matrix $M_\mu$ defined in~(\ref{eq4-2-2}) is positive semidefinite, then, 
$$
\mathrm{sup}\left\{\rho\left(\sum_{k=1}^{N}{A_k\,\mathrm{e}^{j\theta_k}}\right),\; \theta_k\in[0,2\pi]\right\}<1.
$$
\end{cor}
\begin{pf}
Let $\theta_1,\ldots,\theta_N$ be fixed arbitrary reals. Premultiplying and postmultiplying the matrix $M_\mu$, respectively, by
$$
\begin{bmatrix}\mathrm{e}^{-j\theta_1}I_n & \cdots & \mathrm{e}^{-j\theta_N}I_n\end{bmatrix}\; \mathrm{and} \;
\begin{bmatrix}\mathrm{e}^{j\theta_1}I_n \\ \vdots \\ \mathrm{e}^{j\theta_N}I_n\end{bmatrix},
$$
we obtain
\begin{equation}\label{eq4-2-4}
A^\ast P_1 A - P_1 + N \leq 0,
\end{equation}
where $N = \sum_{k=1}^{N}{(\mathrm{e}^{-2\mu r_{k-1}}-\mathrm{e}^{-2\mu r_k})P_k}$, and $A = \sum_{k=1}^{N}{A_k\mathrm{e}^{j\theta_k}}$. Noting that $N$ is a positive definite matrix,~(\ref{eq4-2-4}) implies that
$$
\left(\sum_{k=1}^{N}{A_k\mathrm{e}^{j\theta_k}}\right)^\ast P_1 \left(\sum_{k=1}^{N}{A_k\mathrm{e}^{j\theta_k}}\right) - P_1  < 0,
$$ 
where $P_1$ is symmetric positive definite. Consequently, $\rho\left(\sum_{k=1}^{N}{A_k\,\mathrm{e}^{j\theta_k}}\right)<1$. This inequality being fulfilled for any constants $\theta_1,\ldots,\theta_N$, the result follows.
\hfill{$\Box$}
\end{pf}
It is noted that the converse of this result is false, in general.
\begin{exmp}\label{ex-3}
Take the scalar system
\begin{equation}\label{eq-ex3}
x(t)=0.2 \,x(t-1) -0.05 \,x(t-\sqrt{2}) -0.5\, x(t-2 \pi),
\end{equation}
with initial condition $\varphi(t)= 2 \sin (t) + 1$, for $t \in[-2\pi,0[ $. Its simulation is reported in Fig.~\ref{fig-3}. The conditions of Theorem~\ref{th4-2-2} are fulfilled, with $\mu=0.0609$, $ P_1=  56.8756$, $P_2= 44.7477$ and $P_3=41.6480$. It follows that this system is $L_2$-exponentially stable, with
$$
\|x_t(\varphi)\|_{L_2} \leq 3.7881 \cdot \|\varphi\|_c\cdot \mathrm{e}^{-0.0609 t},\;\text{for}\;t\geq 0,
$$
and $\|\varphi\|_c=3$.
\begin{figure}[h!]\label{fig-3}
\center
\psfrag{y}[c]{\scriptsize$x(t,\varphi)$}
\psfrag{t}{\scriptsize$t$}
\includegraphics[width=9cm]{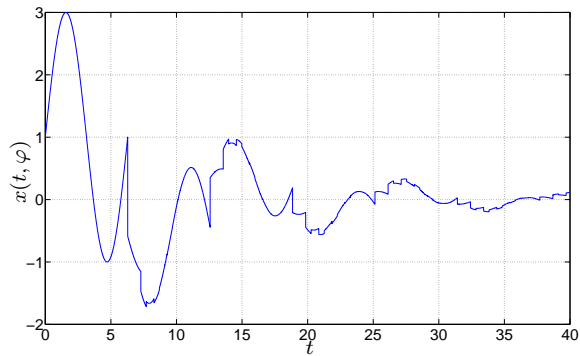}
\caption{Free response $x(t,\varphi)$ of~(\ref{eq-ex3}).}\label{fig-3}
\end{figure}
\end{exmp}

Theorems~\ref{th4-2-1} and~\ref{th4-2-2} are concerned with $L_2$-exponential stability of the solution $x(t,\varphi)$ for~(\ref{eq4-2-1}). Of course, these conditions lead to stability in the delays. It is of interest to know wether exponential stability holds under conditions of Theorem~\ref{th4-2-2}. The positive answer is given in the following corollary, and its proof is reported in Appendix~\ref{appA}. 

\begin{cor}\label{th4-2-4}
Under the conditions of Theorem~\ref{th4-2-2}, exponential stability for $x(t,\varphi)$ holds, that is, for any arbitrarily small $\epsilon\in]0,\mu[$, there exist $\kappa\geq 0$ such that
$$
\|x(t,\varphi)\|\leq \kappa\,\|\varphi\|_c\, \mathrm{e}^{-(\mu-\epsilon) t}\, , \; t\geq 0.
$$
\end{cor}
\begin{pf}
See Appendix~\ref{appA}.
\hfill{$\Box$}
\end{pf}
Few comments on this result can be noticed. While Lyapunov-Krasovskii approach allows to conclude directly on $L_2$-exponential stability, some complementary results issued from the spectral approach were used to conclude on exponential stability. This fact is certainly related to the definition of our Lyapunov-Krasovskii functional, and the absence in~(\ref{eq4-2-1}) of differentiation operator, for which perhaps a more suitable choice should be to analyze variation of Lyapunov-Krasovskii functional and not a differential along the trajectories. It should be also noted that the analysis of discontinuities was not required in the proof of Corollary~\ref{th4-2-4}.   

\section{Robustness under parametric uncertainties}\label{section5}
Until now, stability analysis was only concerned with an exactly known system. It is of interest to give sufficient conditions for stability when norm-bounded parametric uncertainties appear in the model. This section gives a positive answer to such an analysis.

\subsection{Single-delay case }\label{sec5-1}

Consider the uncertain system of the form
\begin{equation}\label{eq5-1-1}
x(t)=(A+ \Delta_A ) x(t-r)
\end{equation}
with some norm-bounded uncertainty matrix $\Delta_A$ such that $ |||\Delta_A|||\leq \delta$,  for some known $\delta\geq 0$, and for $ |||\cdot|||$ the Euclidean subordinated matrix
norm defined by 
\begin{equation}
|||A|||= \underset{||u||=1}{\sup}  ||Au|| = \sqrt{\lambda_{\mathrm{max}}(A^\ast A)}.
\end{equation}
We obtain the following result.

\begin{thm}\label{th5-1-1}
The system~(\ref{eq5-1-1}) is exponentially stable if there exist a symmetric positive definite real matrix $ P $ and $\mu>0$ such that 
\begin{equation} \label{eq5-1-2}
A^T P A - \mathrm{e}^{-2 \mu r} P+\lambda_{\mathrm{max}}(P)(\delta + 2 |||A|||)\delta \leq 0.
\end{equation}
The real $\mu$ is a lower bound for the decay rate of $x(t,\varphi)$, solution of~(\ref{eq5-1-1}).
\end{thm}
\begin{pf}
Assume that~(\ref{eq5-1-2}) holds. Consider the Lyapunov-Krasovskii functional 
$$
v_{\mu}(x_t(\varphi))=\int_{t-r}^{t}{\mathrm{e}^{-2\mu(t-\theta)} x^T(\theta) P x(\theta) \, \mathrm{d}\theta}.
$$
Its time
derivative along the trajectories of~(\ref{eq5-1-1}) is
$$
\dot{v}_{\mu}(x_t(\varphi)) = -2\mu \, v_{\mu}(x_t(\varphi))-x^{T}(t-r) M_{\Delta}  x(t-r)
$$
where 
$$
M_{\Delta}=\mathrm{e}^{-2\mu r}P-A^{T}P A - A^{T}P \Delta_A
 - \Delta_A^{T}P A - \Delta_A^{T} P \Delta_A.
$$
For any element $u\in \Ree^{n}$,
\begin{align}
u^{T} A^{T} P \Delta_A u &\leq ||Au|| \; ||P \Delta_A u || \nonumber \\
&\leq ||u|| \; |||A||| \; |||P||| \; ||\Delta_A u || \nonumber\\
&\leq  \lambda_{\mathrm{max}}(P) \, |||A|||  \; |||\Delta_A||| \; u^Tu.  
\end{align}
Repeating this argument for the last three terms in $M_{\Delta}$, we see that~(\ref{eq5-1-2}) is a sufficient condition to ensure that $M_\Delta\geq 0$. The assertion that $\mu$ is a lower bound for the decay rate of the solution follows from the pf of Theorem~\ref{th4-1-1}.
\hfill{$\Box$}
\end{pf}


%
\begin{exmp}
Take the uncertain system in Example~\ref{ex-2}
\begin{equation}
x(t)= (A+ \Delta_{A})  x(t-\pi),
\end{equation}
where $|||\Delta_A |||\leq \delta =  0.01 $ and $A$ is given in~(\ref{eq-ex2}). A solution of~(\ref{eq5-1-2}) is $\mu=0.2354$ and 
$$P= \begin{pmatrix}
  4.5412  & -2.5013\\
   -2.5013   &  3.5768
\end{pmatrix}. $$
The system is exponentially stable, i.e., 
$$||x(t,\varphi)||\leq   2.0905 \,||\varphi ||_c \,\mathrm{e}^{-0.2354\, t},\; t\geq 0. $$
\end{exmp}

\subsection{Multi-delays case}\label{section5-2}
For the general case, the uncertain system is defined by
\begin{equation}\label{eq5-2-1}
x(t)= \sum_{k=1}^{N}   (A_k+ \Delta_{A_{k}} ) x(t-r_k)
\end{equation}
with $ |||\Delta_{A_{k}}|||\leq \delta_{k}$ and $\delta_{k} \geq 0 $, for $k=1,\ldots,N$. Then, we use this intermediate result.

\newdefinition{lem}{Lemma}
\begin{lem}\label{th5-2-1}
Let $ P $ be a $n\times n$ symmetric positive definite real matrix. Then, for any real vectors $ u $ and $ v $ in $\Ree^n$  and any $n\times n$ real matrices $A$ and $B$, the following inequality holds
\begin{align}
u^{T}A^{T} P B v &\leq \frac{1}{2} \lambda_{\max}(P) \, |||A|||  \,|||B||| (u^Tu+v^Tv). \nonumber
\end{align}
\end{lem}
\begin{pf} 
For any real vectors $u$ and $v$,
\begin{align}
u^T A^T P B v &\leq  ||Au|| \, ||PBv ||,\nonumber \\
&\leq ||u|| \, |||A |||  \, |||P |||  \, |||B||| \, ||v||,\nonumber \\
& = \lambda_{\max}(P) ||u|| \, |||A ||| \, |||B||| \, ||v||.\nonumber
\end{align}
Hence the inequality $\|u\|\|v\| \leq \frac{1}{2} \|u\|^2+  \frac{1}{2} \|v\|^2$ leads to the desired inequality.
\hfill{$\Box$}
\end{pf}

\begin{thm}\label{th5-2-2}
The system~(\ref{eq5-2-1}) is $L_2$-exponentially stable if there exist symmetric positive definite real matrices $P_k$, for $k=1,\ldots,N$, and $\mu>0$ such that  
\begin{equation} \label{eq5-2-2}
-M_{\mu}+\lambda_{\mathrm{max}}(P_{1})  Q_{\Delta}
 \leq 0,
\end{equation}
where $M_\mu$ is given in~(\ref{eq4-2-2}), 
$$
Q_{\Delta}=\mathrm{block\,diag}\{Q_{\Delta_1},\ldots,Q_{\Delta_N}\}, \;\text{with}
$$  
$$
Q_{\Delta_j} =  \sum_{p=1}^{N}{\left(|||A_{j}||| \delta_p+  \delta_j |||A_p ||| + \delta_p \delta_j\right)} \cdot I_n
$$
$n\times n$ matrices, for $j=1,\ldots,N$. The real $\mu$ is a lower bound for the decay rate of the solution $x(t,\varphi)$ of~(\ref{eq5-2-1}). 
\end{thm}
\begin{pf}
Assume that~(\ref{eq5-2-2}) holds. From the Lyapunov-Krasovskii functional
$$
v_\mu(x_t(\varphi))=\sum_{k=1}^{N}\int_{t-r_k}^{t-r_{k-1}} \mathrm{e}^{-2\mu(t-\theta)} x^{T}(\theta) P_k x(\theta) \,\mathrm{d}\theta,
$$
we see that its time derivative along the trajectories of~(\ref{eq5-2-1}) is
$$
\dot{v}_{\mu}(x_t(\varphi)) = -2\mu \,v_{\mu}(x_t(\varphi)) + \psi^T(t) \tilde{M}\psi(t),
$$
where $\tilde{M}= -M_{\mu} + Q $ and 
$$
\psi^T(t)=\begin{bmatrix}
x^T(t-r_1)& \cdots &  x^T(t-r_N)
\end{bmatrix}. 
$$
Denoting $Q_{ij}$ the $n\times n$ entry block in position $(i,j)$ of $Q$, we have 
$$
Q_{ij}=\Delta_{A_i}^T P_1(A_j+\Delta A_j) + A_i^T P_1 \Delta_{A_j},
$$ 
for $i,j=1,\ldots,N$. Applying Lemma~\ref{th5-2-1} for each matrix block $Q_{ij}$, for $i,j=1,\ldots,N$, we see that~(\ref{eq5-2-2}) is a sufficient condition to ensure that $\tilde{M}\leq 0$. The fact that $\mu$ is a lower bound for the decay rate of the solution $x(t,\varphi)$ comes from Theorem~\ref{th4-2-2}.\hfill{$\Box$}
\end{pf}
One can remark, from Corollary~\ref{th4-2-4}, that a similar result holds for exponential stability under uncertainties in parameters. 
\begin{exmp}
Consider the uncertain plant of Example~\ref{ex-3} 
\begin{eqnarray*}\label{param_3}
x(t) & = & (0.2+\Delta_{A_{1}} ) x(t-1) + (- 0.05+\Delta_{A_{2}}) x(t-\sqrt{2})\\ & & + ( -0.5+\Delta_{A_{3}}) x(t-2\pi),
\end{eqnarray*}
with $|||\Delta_{A_{1}}|||\leq 0.01  $, $|||\Delta_{A_{2}}|||\leq 0.03 $ and $|||\Delta_{A_{3}}|||\leq 0.1$. A solution of~(\ref{eq5-2-2}) is $P_1= 16.6281$, $P_2=13.1068$, $P_3=11.7608$ and $\mu = 0.0244$. $L_{2}$-exponential (as well as exponential) stability is then ensured for this family of uncertain systems, with
$$ ||x_t (\varphi) ||_{L_2} \leq  3.0276 \; ||\varphi ||_c \; \mathrm{e}^{-\mu t}, \hspace{0.3 cm}\forall t \geq 0.$$

\end{exmp}

%
\section{Robustness under time-varying delays}\label{section6}
From the previous approach for the characterization of stability, robustness for time-varying delays can be investigated, leading to sufficient conditions of stability. These conditions give a positive answer to the open questions outlined in~\cite{avellarmarconato}.

\subsection{Single delay case}\label{sec6-1}

Consider the system with a single time-varying delay 
\begin{equation}\label{eq6-1-1}
x(t)=Ax(t-r(t)),\; t\geq 0
\end{equation}
where $r(t)=r_0+ \delta_{r}(t)$ with $r_0>0$ some known constant delay. We assume that $\delta_{r}(t)$ is a continuous and differentiable function and 
$-r_0 < \delta_r (t)\leq \delta$, for $\delta\geq 0 $, and $\dot{\delta}_r(t)\leq \delta_1 <1$, for $\delta_1\in \Ree$, to ensure causality of~(\ref{eq6-1-1}).

%

If we want to characterize the degradation of the exponential decay rate when the delay is time-varying, with some unknown time-dependent part in the delay, we obtain this result.

\begin{thm}\label{th6-1-1}
Assume that 
$$z(t)=A \,z(t-r_0) $$
is exponentially stable, that is there exist a symmetric positive definite real matrix $P$ and $\mu>0$ such that
\begin{equation}\label{eq6-1-2}
-M_{\mu}=A^TPA-\mathrm{e}^{-2\mu r_0}P \leq 0. 
\end{equation}
If 
\begin{equation}\label{eq6-1-3}
\dot{\delta}_r(t) \leq \delta_1 < 1-\beta-\mathrm{e}^{-2\mu r_0},
\end{equation}
where $\beta=\frac{\lambda_{\mathrm{max}}(-M_\mu)}{\lambda_{\mathrm{max}}(P)}$, then the system~(\ref{eq6-1-1}) is exponentially stable, that is there exists $\gamma>0$, such that for any arbitrarily small $\varepsilon>0$, 
\begin{equation}\label{eq6-1-4}
\gamma_{\mathrm{max}}-\frac{\varepsilon}{2(r_0+\delta)} \leq \gamma \leq \gamma_{\mathrm{max}},
\end{equation}
where $\gamma_{\mathrm{max}} = -\frac{1}{2(r_0+\delta)}\mathrm{ln}\left(\frac{\beta+\mathrm{e}^{-2\mu r_0}}{1-\delta_1}\right)$,
and
$$
\|x(t,\varphi)\|\leq \sqrt{\frac{\lambda_{\mathrm{max}}(P)}{\lambda_{\mathrm{min}}(P)}} \|\varphi\|_c \, \mathrm{e}^{-\gamma t},\; t\geq 0.
$$

\end{thm}
\begin{pf}
From assumptions, let us define the Lyapunov-Krasovskii functional for~(\ref{eq6-1-1}) 
$$
v_{\gamma}(x_t(\varphi))=\int_{t-r(t)}^{t} { \mathrm{e}^{-2 \gamma (t-\theta)}} x^T(\theta) P x(\theta)\, \mathrm{d}\theta,
$$
where $P$ satisfies~(\ref{eq6-1-2}). The time derivative of $v_\gamma(x_t(\varphi))$ along the trajectories of~(\ref{eq6-1-1}) is
$$
\frac{\mathrm{d}}{\mathrm{d}t}v_\gamma = -2\gamma v_\gamma + x^T(t-r(t))N_\gamma(t) x(t-r(t)),
$$
where
$$
N_\gamma(t) = -M_\mu +\left[\mathrm{e}^{-2\mu r_0}- (1-\dot{\delta}_r(t))\mathrm{e}^{-2 \gamma r(t)}\right]P.
$$
From~(\ref{eq6-1-2}) and~(\ref{eq6-1-3}), it is noticed that 
$$
0 \leq \frac{\beta+\mathrm{e}^{-2\mu r_0}}{1-\delta_1}<1
$$
holds. 
Hence, for any $\varepsilon>0$ such that
$$
\varepsilon\leq-\mathrm{ln}\left(\frac{\beta+\mathrm{e}^{-2\mu r_0}}{1-\delta_1}\right)-\mathrm{max}\{0,-\mathrm{ln}\left(\frac{\mathrm{e}^{-2\mu r_0}-\beta}{1-\delta_1}\right)\},
$$
there exists $\gamma>0$ such that
$$
\gamma_{\mathrm{max}}-\frac{\varepsilon}{2(r_0+\delta)} \leq \gamma \leq \gamma_{\mathrm{max}},
$$
where $\gamma_{\mathrm{max}} = -\frac{1}{2(r_0+\delta)}\mathrm{ln}\left(\frac{\beta+\mathrm{e}^{-2\mu r_0}}{1-\delta_1}\right)$. It follows that for such a $\gamma$,
$$
\eta = \beta+|\mathrm{e}^{-2\mu r_0} -(1-\delta_1)\mathrm{e}^{-2\gamma (r_0+\delta)} |\leq 0.
$$ 
We conclude that
$$
N_\gamma(t) \leq \eta \,\lambda_{\mathrm{max}}(P)\cdot I_n \leq 0,\;\forall t\geq 0,
$$
so that
$$
\frac{\mathrm{d}}{\mathrm{d}t}v_\gamma(x_t(\varphi)) +2\gamma \,v_\gamma(x_t(\varphi)) \leq 0,\;t\geq 0.
$$
Defining similar lower and upper bounds for $v_\gamma(x_t(\varphi))$ than those used in the pf of Theorem~\ref{th4-2-2}, we conclude on $L_2$-exponential stability for $x(t,\varphi)$ from Theorem~\ref{th4-2-1}. We show next that exponential stability holds. For this, remark that
\begin{eqnarray}
N_\gamma(t) & = & A^TPA-(1-\dot{\delta}_r(t))\mathrm{e}^{-2 \gamma r(t)}P\nonumber\\
& \leq & A^TPA-(1-\delta_1)\mathrm{e}^{-2 \gamma r(t)}P \, \leq \, 0. \label{eq6-1-5}
\end{eqnarray}
Note that by construction $0<1-\beta-\mathrm{e}^{-2\mu r_0}\leq 1$. It is then always possible to take $\delta_1$ in~(\ref{eq6-1-3}) such that $0\leq \delta_1<1$. Then, premultiplying and postmultiplying the inequality~(\ref{eq6-1-5}) by $x^T(t-r(t))$ and $x(t-r(t))$, respectively, leads to 
$$
x^T(t)Px(t) \leq \mathrm{e}^{-2 \gamma r(t)} x^T(t-r(t))Px(t-r(t)),\;t\geq0.
$$
Iterating such inequality as in the pf of Theorem~{\ref{th4-1-1}}, we conclude that
$$
\|x(t,\varphi)\|\leq \sqrt{\frac{\lambda_{\mathrm{max}}(P)}{\lambda_{\mathrm{min}}(P)}} \|\varphi\|_c \, \mathrm{e}^{-\gamma t},\; t\geq 0,
$$
which proves exponential stability.
\hfill{$\Box$}
\end{pf}

\begin{exmp}
Take the plant  
\begin{equation}\label{rr1}
x(t)= \begin{pmatrix}
\frac{1}{2} & -\frac{3}{10}\\
\frac{7}{20} & 0
\end{pmatrix} x(t-r(t)),
\end{equation}
with the initial condition $\varphi(t)=\begin{bmatrix}
\sin(3t)\\
\cos(3t)
\end{bmatrix}$, for $t\in [-\pi-\frac{1}{2},0[$, and $r(t)=r_0+\delta_r (t)$ for $r_0=\pi$, $\delta_r (t)=\frac{1}{2}\sin(\frac{t}{2})$. The system~(\ref{rr1}) is exponentially stable, i.e., 
$$||x(t,\varphi)||\leq  2.7951 \,||\varphi ||_c \,\mathrm{e}^{-\gamma t} $$
where $||\varphi ||_c=2$ and, for any arbitrarily small $\epsilon>0$, 
$$\gamma_{\max} - \frac{\epsilon}{2(r_0+\delta)} \leq  \gamma \leq \gamma_{\max} $$
with $\delta=\frac{1}{2}$ and $\gamma_{\max}= 0.2697$. A simulation of the free response for~(\ref{rr1}) is plotted in Fig.~\ref{fig-5}.\\
\begin{figure}[h!]\label{fig-5}
\center
\psfrag{y}[c]{\scriptsize$x(t,\varphi)$}
\psfrag{t}{\scriptsize$t$}
\includegraphics[width=9cm]{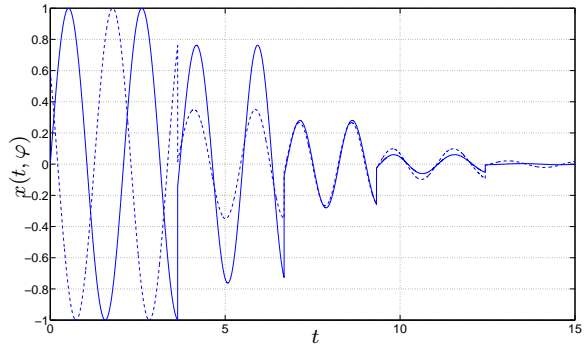}
\caption{Free response $x(t,\varphi)$ of~(\ref{rr1}).}\label{fig-5}
\end{figure}
\end{exmp}

Theorem 18 generalizes the stability result obtained in~\cite{damak2013}, where asymptotic stability in presence of time-varying delays was characterized. Indeed, for asymptotic stability characterization,~(\ref{eq6-1-3}) leads to the sufficient condition (taking $\mu=0$ and $M_\mu$ in~(\ref{eq6-1-2}) positive definite) 
$$
\dot{\delta}_r(t) < \delta_{\mathrm{max}} = \min\left\{1,\frac{\lambda_{\mathrm{min}}(M_\mu)}{\lambda_{\mathrm{max}}(P)}\right\}\,,\, \forall t\geq 0,
$$
which is precisely the condition appearing in~\cite{damak2013} (Theorem 13).

\subsection{Multi-delays case}\label{section6-2}

For the general case of time-varying delays, consider 
\begin{equation}\label{eq6-2-1}
x(t)=\sum_{k=1}^{N} A_{k} x(t-r_{k}(t)),\hspace{0.4 cm} t\geq 0
\end{equation}
where for $k=1,\ldots,N$,
$$
r_{k}(t)=r_{0_k}+ \delta_{r_k} (t) , \hspace{0.5 cm} t \geq 0,
$$
and $\delta_{r_k}(t)$ are
bounded continuous differentiable functions satisfying
$$ 
-r_{0_k}<\delta_{r_k}(t) \leq \delta_{k}, 
$$
with $\delta_k \in \Ree$, and $\dot{\delta}_{r_k}(t)\leq \delta_{1_k}<1$.

In the following result, we take by convention $P_{N+1}=0$.


\begin{thm}
For the system
$$
z(t)=\sum_{k=1}^{N}A_k z(t-r_{0_k}), 
$$
assume that there exist symmetric positive definite real matrices $P_k$, $k=1,\ldots,N$, and $\mu>0$ such that $M_\mu$ in~(\ref{eq4-2-2}) is a positive semidefinite matrix. If for any $k=1,\ldots,N$,
\begin{equation}\label{eq6-2-2}
\dot{\delta}_{r_k}(t) \leq \delta_{1_k} < 1 - \beta -\mathrm{e}^{-2\mu r_{0_{k}}},\;\forall t\geq 0,
\end{equation}
where $\beta= \dfrac{\lambda_{\max}(-M_{\mu})}{\max\limits_{k=1,\ldots,N}\{\lambda_{\max}(P_k) - \lambda_{\min}(P_{k+1})\}} $, then~(\ref{eq6-2-1}) is $L_{2}$-exponentially stable, that is there exist $\gamma>0$ and $\alpha\geq0$ such that for any arbitrarily small $\varepsilon>0$, 
\begin{equation}\label{eq6-2-3}
\gamma_{\mathrm{max}}-\frac{\varepsilon}{2\cdot\max\limits_{k}\{r_{0_k}+\delta_k\}} \leq \gamma \leq \gamma_{\mathrm{max}},
\end{equation}
where $\gamma_{\mathrm{max}} = \min\limits_{k}\left\{-\frac{1}{2(r_{0_k}+\delta_k)}\mathrm{ln}\left(\frac{\beta+\mathrm{e}^{-2\mu r_{0_k}}}{1-\delta_{1_k}}\right)\right\}$, and
$$
\|x_t(\varphi)\|_{L_2} \leq \alpha\|\varphi\|_c\,\mathrm{e}^{-\gamma t}. $$ 
\end{thm}

\begin{pf}
Take the Lyapunov-Krasovskii functional
$$
v_\gamma(x_t(\varphi))=\sum_{k=1}^{N} \int_{t-r_{k}(t)}^{t-r_{k-1}(t)} { \mathrm{e}^{-2 \gamma (t-\theta)}} x^T(\theta) P_k x(\theta)\, \mathrm{d}\theta, 
$$
where $P_{k}$, for $k=1,\ldots ,N$, and $\mu >0$ are the solutions of~(\ref{eq4-2-2}). 
Its time derivative along the trajectories of~(\ref{eq6-2-1}) is given by
$$\dot{v}_\gamma(x_t(\varphi)) = -2\gamma \,v_\gamma(x_t(\varphi))+\psi^T(t)N_\gamma(t)\psi(t),$$
with
$$
\psi(t) = \begin{bmatrix}
x^T(t-r_1(t)) & \cdots & x^T(t-r_N(t)) 
\end{bmatrix}^T.
$$
The matrix $N_\gamma$ is given by
$$
N_\gamma(t) = -M_\mu + Q(t),
$$
where $M_\mu$ is in~(\ref{eq4-2-2}), and, for $P_{N+1}=0$,
$$
Q(t) = \mathrm{block\,diag}\{(\mathrm{e}^{-2\mu r_{0_k}} - (1-\dot{\delta}_{r_k}) \mathrm{e}^{-2\gamma r_{k}} )(P_{k}-P_{k+1})\}.
$$
It is straightforward to verify that, for any $t\geq 0$,
$$
N_\gamma(t)  \leq  (\beta+\chi)\cdot\max\limits_{k}\{\lambda_{\max}( P_{ k}) -\lambda_{\min} (P_{k+1})\}\cdot I_{nN}
$$
where $$
\chi = \max\limits_{k=1,\ldots,N}\left| \, \mathrm{e}^{-2\mu r_{k_{0}}}-(1-\delta_{1k}) \mathrm{e}^{-2\gamma (r_{k_{0}} +\delta_k)} \, \right|.
$$
Using the fact that $\lambda_{\max}( P_{ N})>0$ and $\lambda_{\min} (P_{N+1})=0$, we conclude that $\max\limits_{k}\{\lambda_{\max}( P_{ k}) -\lambda_{\min} (P_{k+1})\}>0$. Hence a sufficient condition to obtain $N_\gamma(t)\leq0$, for any $t\geq$, is $\beta+\chi\leq 0$, or equivalently,
\begin{equation}\label{eq6-2-4}
\left| \, \mathrm{e}^{-2\mu r_{k_{0}}}-(1-\delta_{1k}) \mathrm{e}^{-2\gamma (r_{k_{0}} +\delta_k)} \, \right| \leq -\beta,
\end{equation}
for $k=1,\ldots,N$. We show below that if~(\ref{eq6-2-2}) is satisfied, then~(\ref{eq6-2-4}) holds.\\
For this, from~(\ref{eq4-2-2}), we remark that, for any $k=1,\ldots,N$,    
\begin{eqnarray*}
0 & \leq &  \lambda_{\max}(-M_\mu)+\\
&  &  +\mathrm{e}^{-2\mu r_{0_k}}\max\limits_{k=1,\ldots,N}\{\lambda_{\max}( P_{ k}) -\lambda_{\min} (P_{k+1})\}.
\end{eqnarray*}
Then
$$
\beta+\mathrm{e}^{-2\mu r_{0_k}}\geq 0,\;\mathrm{for}\;k=1,\ldots,N.
$$
Consequently, $\delta_{1_k}<1$ and $0\leq \frac{\beta+\mathrm{e}^{-2\mu r_{0_k}}}{1-\delta_{1_k}}< 1$, for $k=1,\ldots,N$. It follows that $\gamma_{\max}>0$. For any $\varepsilon>0$ such that, for $k=1,\ldots,N$,
$$
\varepsilon\leq-\mathrm{ln}\left(\frac{\beta+\mathrm{e}^{-2\mu r_{0_k}}}{1-\delta_{1_k}}\right)-\mathrm{max}\{0,-\mathrm{ln}\left(\frac{\mathrm{e}^{-2\mu r_{0_k}}-\beta}{1-\delta_{1_k}}\right)\},
$$
there exists $\gamma>0$ such that
$$
\gamma_{\mathrm{max}}-\frac{\varepsilon}{2\cdot\max\limits_{k}\{r_{0_k}+\delta_k\}} \leq \gamma \leq \gamma_{\mathrm{max}}.
$$
For such a $\gamma>0$, the inequality (\ref{eq6-2-4}) holds. We then conclude that $\dot{v}_\gamma+2\gamma\, v_\gamma \leq 0$, for any $t\geq 0$.
Similarly to Theorem~\ref{th4-2-2}, we have
$$
\alpha_1\|x_t(\varphi)\|_{L_2}^2 \leq v_\gamma(x_t(\varphi))
$$
with $\alpha_1=\min\limits_{k}\{\mathrm{e}^{-2\gamma(r_{0_k}+\delta_k)}\lambda_{\mathrm{min}}(P_k)\}>0$, and
$$
v_\gamma(\varphi) \leq \alpha_2 \|\varphi\|_c^2,
$$
with $\alpha_2 = (r_{0_N}+\delta_N)\cdot\max\limits_{k}\lambda_{\mathrm{max}}(P_k)$. $L_2$-exponential stability of $x(t,\varphi)$ follows from~Theorem~\ref{th4-2-1}, with $\alpha=\sqrt{\frac{\alpha_2}{\alpha_1}}$.\hfill{$\Box$}
\end{pf}

\begin{exmp}\label{ex21}
Consider the plant
\begin{equation*}\label{eq-ex6}
x(t)= 0.2\, x(t-r_1(t))-0.05\, x(t-r_2(t))-0.5\, x(t-r_3(t)),
\end{equation*}
with an initial condition $\varphi(t)=2 \sin(t)+1$, for $t\in [-(2\pi+1),0[$, $r_1(t)= 1+ \delta_{r_{1}}(t)$, $r_2(t)= \sqrt{2}+ \delta_{r_{2}}(t)$, $r_3(t)= 2\pi + \delta_{r_{3}}(t)$, $\delta_{r_{1}}(t)=\frac{1}{2}\,\sin(\frac{t}{5})$, $\delta_{r_{2}}(t)=0.15\,\sin(t)$ and $\delta_{r_{3}}(t)=\sin(0.4\, t)$. Its simulation is reported in Fig.~\ref{fig-6}. Taking $\delta_{1_1}=0.1$, $\delta_{1_2}=0.15$, $\delta_{1_3}=0.4$, a solution of~(\ref{eq6-2-2}) is obtained with
$$
\beta = -6.01\cdot 10^{-5},\, \gamma_{\mathrm{max}}=0.0031.
$$
It follows that this system is $L_2$-exponentially stable, with
$$
\|x_t(\varphi)\|_{L_2} \leq   4.9125  \, ||\varphi ||_c\, \mathrm{e}^{-\gamma t},\;\text{for}\;t\geq 0,
$$
with $\|\varphi\|_c=2$ and $\gamma$ obtained from~(\ref{eq6-2-3}).
\begin{figure}[h!]\label{fig-6}
\center
\psfrag{y}[c]{\scriptsize$x(t,\varphi)$}
\psfrag{t}{\scriptsize$t$}
\includegraphics[width=9cm]{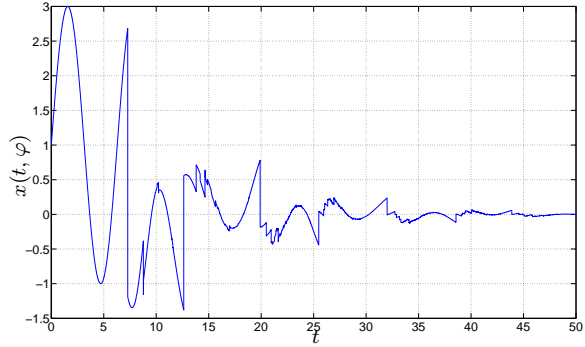}
\caption{Free response $x(t,\varphi)$ in Example~\ref{ex21}.}\label{fig-6}
\end{figure}
\end{exmp}

\appendix
\section{Appendix}\label{appA}    

\textbf{proof of Corollary~\ref{th4-2-4} : } The proof is divided in two steps. In the first step, we show that under the conditions of Theorem~\ref{th4-2-2}, the solution $x(t,\varphi)$ is bounded for any $t\geq 0$. 

Consider the distribution
$$
g(t) = \delta(t) - \sum_{k=1}^{N}{A_k\delta(t-r_k)},
$$
where $\delta(t)$ stands for the Dirac distribution. This distribution lies in the Banach algebra $\ell_1^{n\times n}$, which is the algebra of elements in the form $f(t)=\sum_{i\geq 0}{f_i \delta(t-t_i)}$, where $0=t_0<t_1<\ldots$, $f_i\in\Cee^{n\times n}$, $\sum_{i\geq 0}{\|f_i\|_1}<\infty$. The notation $\|\cdot\|_1$ stands for the 1-induced matrix norm.\\
Assume that the conditions of Theorem~\ref{th4-2-2} are fulfilled. From~Corollary~\ref{th4-2-3}, we know that
$$
\mathrm{sup}\left\{\rho\left(\sum_{k=1}^{N}{A_k\,\mathrm{e}^{j\theta_k}}\right),\; \theta_k\in[0,2\pi]\right\}<1.
$$
We show below that under such an assumption, $g$ is a unit of $\ell_1^{n\times n}$, that is $g$ is invertible and its inverse is in $\ell_1^{n\times n}$.\medskip\\
From~\cite{desoervidyasagar}, we know that $g\in\ell_1^{n\times n}$ is invertible over $\ell_1^{n\times n}$ if and only if
$$
\mathop{\mathrm{inf}}_{\mathrm{Re\,s\geq 0}}|\mathrm{det}\,\hat{g}(s)|>0,
$$
where $\hat{g}(s)$ is the Laplace transform of $g(t)$. By contradiction, two cases may arise. First, assume that there exists $\lambda_0=x_0+jy_0$ with $x_0\geq 0$, such that $|\mathrm{det}\,\hat{g}(\lambda_0)|=0$, that is
$$
|\mathrm{det}(I_n-\sum_{k=1}^{N}{A_k\mathrm{e}^{-x_0r_k}\mathrm{e}^{-jy_0r_k}})|=0.
$$ 
This last equality implies that
$$
\rho(\sum_{k=1}^{N}{A_k\mathrm{e}^{-x_0r_k}\mathrm{e}^{-jy_0r_k}})\geq 1.
$$
From Corollary~\ref{th4-2-3}, this implies that there does not exist symmetric positive definite matrices $Q_k$ and $\eta>0$ such that $M_{\eta}$ in~(\ref{eqAPP-1}) is positive semidefinite. 
\begin{figure*}[ht]
\begin{equation}\label{eqAPP-1}
-M_{\eta} = \begin{bmatrix}
\mathrm{e}^{-2x_0r_1}A^{T}_{1} Q_{1} A_{1} + {\mathrm{e}^{-2 \eta r_{1}}} (Q_2 - Q_{1})     &  \cdots  & \mathrm{e}^{-x_0(r_1+r_N)}A_{1}^{T} Q_{1} A_{N}\\
\mathrm{e}^{-x_0(r_1+r_2)}A^{T}_{2} Q_{1} A_{1}    &   \ddots    & \mathrm{e}^{-x_0(r_2+r_N)}A^{T}_{2} Q_{1} A_{N}  \\
 \vdots & \ddots & \vdots\\
  \vdots   & \ddots   &\vdots\\
\mathrm{e}^{-x_0(r_1+r_N)}A^{T}_{N} Q_{1} A_{1}  & \cdots & \mathrm{e}^{-2x_0r_N}A^{T}_{N} Q_{1} A_{N}-{\mathrm{e}^{-2 \eta r_N}} Q_N 
\end{bmatrix}
\end{equation}
\hrulefill
\end{figure*}
Nevertheless, from assumption, $-M_\mu$ in~(\ref{eq4-2-2}) is negative semidefinite, for some $P_k$, $k=1,\ldots,N$ and $\mu>0$. We also remark that
\begin{equation}\label{eqAPP-2}
-M_\eta = - D_0 M_\mu D_0
\end{equation} 
with $D_0 = \mathrm{diag}(\mathrm{e}^{-x_0r_1}I_n,\ldots,\mathrm{e}^{-x_0r_N}I_n)$ a positive definite matrix,  $\eta = \mu+x_0>0$ and $Q_k=P_k$. Such $M_\eta$ is positive semidefinite, which leads to a contradiction.\medskip\\
The second case that may arise occurs for a sequence of $\lambda_q=x_q+jy_q$ such that $|\mathrm{det}\,\hat{g}(\lambda_q)|=0$ with $x_q\rightarrow 0$. Similarly to the previous case, this implies that
$$
\rho(\sum_{k=1}^{N}{A_k\mathrm{e}^{-x_qr_k}\mathrm{e}^{-jy_qr_k}})\geq 1.
$$
By continuity of the spectral radius, we conclude that there exists $y$ such that
$\rho(\sum_{k=1}^{N}{A_k\mathrm{e}^{-jyr_k}})\geq 1$. This is a contradiction. So $g$  is a unit of $\ell_1^{n\times n}$.\medskip\\
We show next that the solution $x(t,\varphi)$ of~(\ref{eq4-2-1}) is bounded for any $t\geq 0$. For this, apply the Laplace transform to (\ref{eq4-2-1}) to get
$$
(I_n-\sum_{k=1}^{N}{A_k \mathrm{e}^{-sr_k}}) \hat{x}(s) = \hat{\psi}(s),
$$
where $\hat{\psi}(s) = \sum_{k=1}^{N}{A_k \mathrm{e}^{-sr_k}\hat{\varphi}_k(s)}$, and $\hat{\varphi}_k(s) = \int_{-r_k}^{0}{\mathrm{e}^{-su}\varphi(u)\,\mathrm{d}u}$.
In the time-domain, this equation reads 
$$
g(t)\ast x(t) = \psi(t),
$$
where $\psi(t)$ is a function with bounded variation. Since $g$  is a unit of $\ell_1^{n\times n}$, we see that there exists $h\in\ell_1^{n\times n}$ such that
$$
x(t) = h(t)\ast\psi(t),
$$
that is $x(t,\varphi)$ is bounded, for any $t\geq 0$.\medskip\\ 
For any arbitrarily small $\epsilon$ in $]0,\mu[$, define $\gamma=\mu-\epsilon$, and $z(t,\varphi)=\mathrm{e}^{\gamma t}\,x(t,\varphi)$ for $t\geq 0$. Following the same reasoning done in~(\ref{eqAPP-2}), and using the first step of this proof, we know actually that $z(t,\varphi)$ is bounded, so that exponential stability for $x(t,\varphi)$ holds.

\bibliographystyle{model1-num-names}
\bibliography{SDetal_arxiv}

\end{document}